\input amstex

\magnification=1190

\centerline{\bf A characterization of spherical polyhedron surfaces}
\medskip
\centerline{\bf Feng Luo}

%
\medskip
\noindent
\centerline{\it Abstract} 

\medskip
A spherical polyhedron surface is a triangulated surface obtained by isometric gluing of spherical triangles.
 For instance, the boundary of a generic convex polytope in the 3-sphere is a 
spherical polyhedron surface. This paper investigates these surfaces from the point of view of inner angles. 
A rigidity result is obtained. A characterization of spherical polyhedron surfaces in terms of the triangulation and the angle assignment is established.
\medskip
\noindent
\S1. {\bf Introduction}
\medskip
\noindent
1.1.  In an attempt to understand the geometric triangulations of closed 3-manifolds with constant sectional 
curvature metrics, we are led to the study of 
spherical polyhedron surfaces. These are metrics obtained by taking a finite collection of spherical triangles and identifying their edges in pairs by isometries. In particular, they are spherical cone metrics on a surface together with a geometric triangulation.
For instance, the link of a vertex in a 3-dimensional geometric triangulation is a spherical polyhedron surface.
In [Lu1], we have initiated an approach to find constant curvature metrics on triangulated closed 3-manifolds
using dihedral angles as parameters. This leads us to investigate spherical polyhedron surfaces from the inner angle point of view. For a spherical polyhedron surface, its \it edge invariant \rm associates each edge of the triangulation the sum of the two inner angles facing the edge. 
The main result of the paper gives a characterization of the spherical polyhedron metrics in terms of the edge invariant.
To be more precise, we prove that if two spherical polyhedron surfaces with isomorphic triangulations have the same edge invariant, then they are isometric.  We also establish an existence result on spherical polyhedron surfaces when the
edge invariants take values  in $[0, \pi]$. 
Similar results for Delaunay triangulations of  surfaces in the Euclidean or hyperbolic cone metrics
have been  worked out beautifully by Rivin [Ri1] and Leibon [Le].  Our approach follows the strategies in [Ri1], [Le] by using a different energy function.

\medskip
\noindent
1.2. We now set up the frame work. Suppose $S$ is  a closed surface and $T$  is a triangulation of the surface. Here by a
triangulation we mean the following: take a finite collection of triangles and identify their edges in pairs by homeomorphisms.
Let V, E, F be the sets of all vertices, edges and triangles in the triangulation $T$ respectively. 
If $a, b$ are two simplices in the triangulation $T$, we use $a<b$ to denote that $a$ is a face of $b$.
The set of \it corners \rm of  $T$ is
$\{ (e, f) | e \in E, f \in F $ so that $e<f$$\}$ and is denoted by $C(S,T)$. By a \it spherical angle structure \rm on the
 triangulated surface
$(S, T)$ we mean a map $x: C(S,T) \to (0, \pi)$ so that for each $f \in T$ and the three edges $e_1, e_2, e_3$ of $f$,
the numbers $x_i = x(e_i, f)$, $i=1,2,3$, form the inner angles of a spherical triangle.  A \it spherical polyhedron metric  \rm on
the triangulated surface $(S, T)$ is a map $ l: E \to (0, \pi)$ so that for each triangle $f$ and its three edges $e_1, e_2, e_3$,
the three numbers $l_i = l(e_i)$, $i=1,2,3,$ form the edge lengths of a spherical triangle. Evidently, given any
spherical polyhedron metric, there is a natural spherical angle  structure associated to it by measuring its inner angles.  One of the goal
 in the
paper is to
characterize the set of all spherical polyhedron metrics inside the space of all spherical angle  structures. 
To this end, we introduce the notion of  the \it edge invariant \rm $D_x$ of the spherical angle  structure $x$. The
 edge invariant $D_x$ is the map defined on the set of all edges $E$ so that the its value at an edge is the sum of the two
 inner angles facing the edge,
i.e., $D_x(e) = x(e,f) + x(e,f')$ where $f,f' \in F$ and $e<f, e<f'$ (it may occur that $f=f'$). 

\medskip
\noindent
{\bf Theorem 1.1.}  \it Given any triangulated closed surface and a real valued function $D$ defined on the set
 of all edges of the triangulation, there is at most one 
spherical polyhedron metric  having $D$ as the edge invariant. \rm

\medskip

An interesting consequence of theorem 1.1 says that if two convex
spherical polytopes in $S^3$ have the same combinatorial triangulation so that their edge invariants are the same, then
these two polytopes are isometric in $S^3$. 

\medskip

\noindent
{\bf Theorem 1.2.} \it Given any triangulated closed surface and a function $D: E \to (0, \pi)$
 so that there is a spherical angle  structure having $D$ as the edge invariant, then there exists a spherical polyhedron metric  having $D$ as the 
edge invariant function. \rm

\medskip

The existence of spherical angle structures with given edge invariant is a linear programming problem and can be checked
algorithmically. The following theorem has been proved by R. Guo [Gu].

\medskip
\noindent
{\bf Theorem 1.3(Guo [Gu]).} \it  Given any triangulated closed surface and any  function $D: E \to (0, \pi)$,
 there is a spherical angle  structure having $D$ as the edge invariant if and only if for any subset $X$ of triangles
in the triangulation,
$$  \pi |X| < \sum_{ e \in E(X)} D(e), $$
where $E(X)$ is the set of all edges of triangles in $X$ and $|X|$ is the number of triangles in $X$. \rm

\medskip

We remark that a slightly stronger version of theorem 1.2 can also be established for  edge invariants $D(E) \subset (0, \pi]$. See
theorem 2.1. 

\medskip

The space of all spherical polyhedron metrics on $(S, T)$, denoted by $CM(S,T)$  is an open convex polytope of dimension $|E|$,
the number of edges. The space
of all positive functions on the set of all edges $E$ is denoted by $\bold R_{>0}^E$.
The map $\Pi:
CM(S, T) \to \bold R_{>0}^E$ sending a cone metric to its edge invariant is evidently a smooth map between two open cells
of the same dimension. Theorem 1.1 shows that the map is injective (in fact it is a local diffeomorphism). Theorems 1.2 and 1.3 show that the image of
the subset of $CM(S, T)$ with edge invariant $D: E \to (0, \pi)$ under $\Pi$ 
is convex polytope. An interesting question is whether 
 the image of $\Pi$ is  an open convex polyhedron in $\bold R^E_{>0}$?
The situation is a bit similar to Thurston's proof of his circle packing theorem for triangulated surface of negative 
Euler characteristic ([Th]). 

The strategy of proving theorems 1.1 and 1.2 goes as follows. For each spherical triangle, we introduce the
concept of \it capacity \rm of the triangle. The capacity is a strictly convex function defined on the space of all spherical triangles parametrized by the inner angles. We define the capacity of a spherical angle structure to be the sum of the capacities of its triangles.
Then the  capacity defines a strictly convex function on
the space $AS(S,T)$ of all spherical angle  structures on $(S, T)$. Given an edge invariant $D: E \to (0, \infty)$, we consider 
the subset $AS(S, T; D)$ of $AS(S, T)$
consisting of all spherical angle  structures with $D$ as the edge invariant.  We prove that
the critical points of the capacity function restricted to the subspace  $AS(S, T; D)$ 
are exactly  the spherical polyhedron metrics on $(S, T)$. Since a strictly convex function cannot have more than one critical points,
theorem 1.1 follows. For theorem 1.2, we show that the capacity function which has a natural continuous extension to
the compact closure of $AS(S, T; D)$ cannot achieve its minimal points in the boundary. Thus the minimal point of the capacity exists in $AS(S, T; D)$ when $D : E \to (0, \pi)$.

\medskip
\noindent
1.3. The study of geometric structures on triangulated surfaces from variational point of view  has appeared in many works [BS], [Co], [Le], [Ri1] and others. In [Ri1] Rivin studied the Euclidean cone metrics and Leibon [Le] worked out the Delaunay triangulations for hyperbolic surfaces. Results similar to
theorems 1.1, 1.2 and 1.3 were proved for Euclidean and hyperbolic geometric triangulations
 in [Ri1] and [Le]. The approach in this paper follows the work in [Ri1] and [Le]. In  [Ri1] and [Le], the "capacity" of a Euclidean and a hyperbolic triangle was introduced. They are all related to the volume in hyperbolic spaces. In turns out the capacities introduced in [Ri1], [Le] and in our current work
can be summarized in one sentence. Namely, given a spherical, or a Euclidean or a hyperbolic triangle in the Riemann sphere considered as the
infinity of the hyperbolic 3-space, there are three circles bounding the triangle. The capacity of the triangle is essentially (up to multiplication and
addition of constants) the hyperbolic volume of the convex hull of the intersection points of these three circles.
The explicit expressions of the capacities are (3.9) and (3.10). For Euclidean triangles, the ideal hyperbolic convex polytopes
are ideal tetrehedra; for hyperbolic triangles, they are ideal hyperbolic prisms; and for
spherical triangles, they are  ideal hyperbolic
octahedra. In our case, we first discovered the capacity of a spherical triangle through the derivative of the cosine law
and later realized that it is again a hyperbolic volume.
It turns out for spherical triangle, Peter Doyle [Le] defined a different capacity (see (3.10)). Doyle's capacity of a spherical triangle is the volume
of the hyperbolic tetrahedron which is the convex hull of four points consisting of the three vertices of the spherical triangle and the Euclidean center 
in the Poincare model (where the spherical triangles are bounded by great circles).

From this point of view, given a triangulated surface $(S, T)$, there are five linear programming problems and variational problems associated to the surface. The linear programming problems are related to the angle structures and the variational problems are the critical points of the "capacities". To begin, let us introduce some concepts.
An \it angle structure \rm on a triangulated surface $(S, T)$ assigns each corner of $(S, T)$ a number in $(0, \pi)$, called the inner angle.
A \it hyperbolic (or spherical, or Euclidean)
 angle structure \rm is an angle structure so that each triangle with the angle assgnements is hyperbolic (or spherical, or Euclidean).  Euclidean angle structures were  first defined by Rivin in [Ri1] who called it \it locally Euclidean structures. \rm  The basic examples of hyperbolic (or spherical, or Euclidean) angle structures are hyperbolic (or spherical, Euclidean) cone metrics with a geometric triangulation by measuring the inner angles.
Given an angle structure $x: C(S,T) \to \bold R_{>0}$, we define its \it edge invariant, \rm denoted by $D_x: E \to
\bold R_{>0}$ to be the sum of two opposite facing angles and its \it Delaunay invariant \rm $\Cal D_x: E \to \bold R_{>0}
$ to be $\Cal D_x(e) = c+d+f+g -a -b$ where $a,b$ are the two angles facing the edge $e$ and $c,d,f,g$ are the four angles
having $e$ as an edge. An angle structure is called \it Delaunay \rm if its Delaunay invariant $\Cal D_x$ is non-negative.
For Euclidean angle structures, the Delaunay invariant and the edge invariant are related by $ 2D_x + \Cal D_x = 2 \pi$. For a spherical, Euclidean or hyperbolic cone metric with a geometric triangulation, its underlying angle structure is Delaunay if and only if the triangulation satisfies the empty circumcircle property, i.e., the interior of the circumcircle of each triangle does not contain any vertices. 

The five linear programming problems associated to the triangulated surface $(S, T)$ are as follows. Namely, the spaces of all
hyperbolic angle structures with prescribed edge invariant $D$ or Delaunay invariant $\Cal D$,  the spaces of all 
spherical angle structures with prescribed edge invariant $D$ or Delaunay invariant $\Cal D$, and the space of all Euclidean angle structures with prescribed Delaunay invariant $\Cal D$. We denote these five convex polytopes by $AH(S, T; D)$, $AH(S, T; \Cal D)$, $AS(S, T;D)$, $AS(S, T; \Cal D)$ and $AE(S, T; \Cal D)$. In the recent work of R. Guo [Gu], he has found the necessary and sufficient conditions for these spaces to be non-empty. The works of Rivin and Leibon dealt with the spaces $AE(S, T;D)$
and $AH(S, T; \Cal D)$ and used the capacity given by formula (3.10). Our paper addresses the space $AS(S, T;D)$ using capacity (3.9). There remain the  problems on the existence and 
uniqueness  of constant curvature cone metrics in the spaces $AS(S, T; \Cal D)$ and $AH(S, T; D)$. We remark that the
associated energies for these problems have been found. Namely, for the space $AS(S, T; \Cal D)$,  Peter Doyle [Le] associated the capacity function given by (3.10)  and observed that the critical points of the capacity are exactly the spherical cone metrics with geodesic triangulations.
The capacity function for the space $AH(S, T; D)$ is given by (3.9) and it is easy to prove that the critical points of the energy are the hyperbolic cone metrics. However, in both  cases the capacity functions are no longer convex or concave.  It is a very interesting problem to establish the existence of the critical points of the capacity function in these cases. Furthermore, it is also interesting to know if
the critical points are unique in the case of hyperbolic cone metrics in 
$AH(S, T; D)$.

\medskip
In our recent work [Lu1], we proposed a generalization of the above set up for closed triangulated 3-manifolds by introducing the 3-dimensional angle structure and its volume. The link of a vertex in a 3-dimensional
angle structure is a spherical angle structure on the 2-sphere. We tend to think that the Delaunay condition for angle structures in dimension-3 is related to the edge invariant $D$ being in the interval $[0, \pi]$ for surfaces.  This is the motivation of the study in this paper. Another motivation of the study is that a spherical angle structure on surface is a 2-dimensional simple model of the 3-dimensional
project in [Lu1].  Theorems 1.1 and 1.2 give some positive evidences for the 3-dimensional project in [Lu1].  In this  comparision, the resolution of the Milnor conjecture on volume of simplexes in [Lu2] ([Ri2] has a new proof) can be considered as the counterpart of proposition 3.1 in dimension-3.

\medskip






\medskip
\noindent
1.4. The paper is organized as follows. In section 2, we recall some known facts about the derivatives of the cosine laws.
We also introduce the capacity function. Some of the basic properties of the capacity function are established. 
In particular, we prove theorems 1.1 and 1.2 in section 2 assuming two important properties of the capacity function. These two 
properties are established in sections 3 and 4.
 In section 3, we show that the capacity function has a continuous extension to the degenerated spherical triangles by relating it
to the
Lobachevsky function. 
In section 4, we study the behavior of the derivative of the capacity function at the  degenerated spherical triangles. 
\medskip
\noindent
1.4. \it Acknowledgement \rm   I would like to thank the referee for his/her suggestions on improving the exposition of the paper.
The work has been supported in part by the NSF and a research grant from Rutgers University.

\medskip
\noindent
\S2. {\bf  Spherical Triangles and Proofs of Theorems 1.1 and 1.2}

\medskip
\noindent
We prove theorems 1.1 and 1.2 assuming several technical properties on spherical triangles in this section.
For simplicity, we assume that the indices $i,j,k$ are  pairwise distinct in this section.
\medskip
\noindent
2.1.   Given a spherical, Euclidean or hyperbolic triangle with inner angles $x_1, x_2, x_3$, let $y_1, y_2, y_3$ be the edge lengths so that
$y_i$-th edge is facing the angle $x_i$. Let $\lambda=0,-1,1$ be the curvature of the underlying space, i.e., $\lambda=1$ for spherical triangles, $\lambda =-1$ for hyperbolic triangles and $\lambda=0$ for Euclidean triangles. 
The cosine law states that,

$$ \cos (\sqrt{\lambda} y_i) =\frac{ \cos  x_i + \cos  x_j \cos  x_k}{ \sin   x_j \sin   x_k },  \tag 2.1$$
 where $\{i,j,k\}=\{1,2,3\}$.

Furthermore, the partial derivatives of $y_i$ as a function of $x=(x_1, x_2, x_3)$ are given by the following lemma.
\medskip
\noindent
{\bf Lemma 2.1.} \it For any spherical or hyperbolic triangle of inner angles $x_i, x_j, x_k$ and the corresponding edge lengths $y_i, y_j, y_k$, where
$\{i,j,k\}=\{1,2,3\}$, the following hold. 

(a) $\partial y_i/\partial x_i = \sin (x_i) /A_{ijk}$ where $A_{ijk} = \sin( \sqrt{\lambda} y_i) \sin x_j \sin x_k /\sqrt{\lambda} $ satisfies $A_{ijk}=A_{jki}$,

(b) $\partial y_i/\partial x_j =  \partial y_i/\partial x_i \cos y_k$. \rm

\medskip
The proof is a simple exercise in calculus, see for instance [Lu1].   

 The space of all spherical triangles parametrized by its inner angles
$x_1, x_2, x_3$, denoted by  $M_3$, is the open tetrahedron $\{ x=(x_1, x_2, x_3) \in (0, \pi)^3 | x^*_i >0, \sum_{i=1}^3 x_i > \pi$\} where $x^*_i = 1/2( \pi + x_i - x_j -x_k)$, $\{i,j,k\}=\{1,2,3\}$.  To see that these inequalities are necessary, we first note that the sum of inner angles of a spherical triangle is larger than $\pi$. To see
$x_1^* >0$, we note that  if $x_1, x_2, x_3$ are the inner angles of a spherical triangle $A$, then $x_1, \pi-x_2, \pi-x_3$ also form the inner angles of a spherical triangle $B$ so that $A \cup B$ forms a region bounded by two great circles intersecting at an angle $x_1$. It follows that the sum $x_1 + \pi-x_2 + \pi-x_3 > \pi$. This shows $x_1^*>0$ is necessary. It is not difficult to show that these four inequalities are also sufficient.  

\medskip
\noindent
{\bf Corollary 2.2.} \it 
(a) The differential 1-form $w = \sum_{i=1}^3 \ln \tan(y_i/2) dx_i$ is closed in the open set $M_3$.

(b) The function $\theta(x) = \int_{(\pi/2, \pi/2, \pi/2)}^x w$ is well defined on $M_3$ and is strictly convex. 

(c)  The differential 1-form $\tilde w = \sum_{i=1}^3 \ln \tanh(y_i/2) dx_i$ is closed in the  set 
$\Cal H_3 =\{(x_1, x_2, x_3) \in \bold R^3 | x_i >0, \quad x_1+x_2+x_3 < \pi\}$ of all hyperbolic triangles.

(d) The function $\tilde \theta(x) = \int_{(0,0,0)}^x \tilde w$ is a well defined smooth function on $\Cal H_3$. \rm

\medskip
\noindent
{\bf Proof.} To show part (a),  it suffices to prove $\partial ( \ln \tan(y_i/2) )/\partial x_j$ is symmetric in $i,j$. By lemma 2.1,
the partial derivative is found to be 
$$1/\sin (y_i)  \partial y_i/\partial x_j  = \cos(y_k)  [\sin(x_i) / \sin (y_i) ]/A_{ijk}.  \tag 2.2$$
By the sine law, one sees clearly that the partial derivative is  symmetric in
$i,j$.  Note also that 
$$ \partial  (\ln \tan (y_i/2)) /\partial x_i  =  [\sin(x_i) / \sin (y_i) ]/A_{ijk}.  \tag 2.3$$

Since the space $M_3$ is simply connected, we see that the function $\theta(x)$ is well defined on $M_3$.
To show that the function $\theta$ is strictly convex, let us calculate its Hessian matrix $H = [h_{rs}]_{3 \times 3}$.
By definition, we have $h_{rs} = \partial ( \ln \tan(y_r/2) )/\partial x_s$.   By (2.2) and (2.3), we have
$h_{ij} = h_{ii} \cos y_k$ and $h_{11}= h_{22}=h_{33} >0$ by the sine law. Thus the matrix $H$ is a positive
multiplication of the matrix $[a_{rs}]$ where $a_{ij} = \cos y_k$ and $a_{ii} = 1$. For a spherical triangle of edge lengths
$y_1, y_2, y_3$, the matrix $[a_{rs}]$ is always positive definite. Indeed, let $v_1, v_2, v_3$ be the three unit vectors
in the 3-space forming the vertices of the spherical triange, then by definition, $a_{rs}$ is the inner product of $v_r$ with $v_s$.
 Thus
the matrix $[a_{rs}]$ is positive definite since it is the Gram matrix of three independent vectors. 

The verifications of parts (c) and (d) are similar and will be omitted. QED

\medskip
\noindent
2.3. The closure of $M_3$ in $\bold R^3$
is given by $\bar M_3 =\{ x \in [0, \pi]^3 | x^*_i \geq 0, x_1+x_2 + x_3 \geq \pi\}$. In sections 3 and 4, we will establish the following two properties concerning the function $\theta$.  Recall that the Lobachevsky function $\Lambda(t) = - \int_0^t \ln |2 \sin u | du$. The function is continuous
on the real line $\bold R$ and is an odd periodic function of period $\pi$. See Milnor [Mi] for more details.

\medskip
\noindent
{\bf Proposition 3.1.} \it  The capacity function $\theta(x) =
  \int_{(\pi/2, \pi/2, \pi/2)}^x \sum_{i=1}^3 \ln \tan (y_i/2) dx_i$ is given by the following,

$$ \theta (x_1, x_2, x_3) = - \sum_{i=1}^3 \Lambda (x^*_i) - \Lambda(  (\pi +x_1+ x_2 + x_3)/2)  + 4 \Lambda(\pi/4)$$
and  the capacity function $\tilde \theta(x) =\int_{(0,0,0)}^x \sum_{i=1}^3 \ln \tanh(y_i/2) dx_i$ is given by
$$\tilde  \theta (x_1, x_2, x_3) = - \sum_{i=1}^3 \Lambda (x^*_i) - \Lambda((\pi+x_1+ x_2 + x_3)/2). $$
In particular, both $\theta$ and $\tilde \theta$ have continuous extensions to the closure $\bar M_3$ of the moduli space of spherical triangles
$M_3  =\{(x_1, x_2, x_3) \in (0, \pi)^3 | x_1+x_2+x_3 > \pi$ and $x_i^* >0$, $i=1,2,3$\} and the closure of $\{(x_1, x_2, x_3)
\in (0, \pi)^3 | x_1+x_2+x_3 < \pi$\}. 
Geometrically, $16 \Lambda(\pi/4) -4\theta(x_1, x_2, x_3)$ is the volume of the hyperbolic ideal octahedron whose vertices are the intersection points of the three circles bounding the spherical triangle $(x_1, x_2, x_3)$.\rm

\medskip
\noindent
{\bf Proposition 4.1.} \it 
For any point $a \in \bar M_3 - M_3$ and a point $p \in M_3$, let $f(t)$ be the function $\theta ( (1-t) a + t p)$
where $t \in (0,1)$. If $a$ is not one of $(0,0,\pi)$$, (0, \pi, 0),$$ (\pi, 0, 0), $$(\pi, \pi, \pi)$, then
$$\lim_{t \to 0^+} f'(t) = -\infty. $$
If $a \in \{ (0,0,\pi), (0, \pi, 0), (\pi, 0, 0), (\pi, \pi, \pi)\}$, then the limit $\lim_{t \to 0^+} f'(t)$ exists and is a finite number. 
 \rm
\medskip

In the rest of the section, we prove theorems 1.1 and 1.2 assuming propositions 3.1 and 4.1.

\medskip
\noindent
2.4.
Given a spherical triangle of inner angles $x_1, x_2, x_3$, we define its  \it capacity \rm to be $\theta(x_1, x_2, x_3)$ 
where $\theta$ is the function introduced in corollary 2.2.  For a spherical angle  structure, we define its capacity to be
the sum of the capacities of its spherical triangles.
To write down the capacity function explicitly,
let us fix some
notations. First, let us label the set of all corners in $(S, T)$ by integers $\{1, ..., n\}$. If three corners
 labeled by $a,b,c$ are of the form $(e_1, f), (e_2, f), (e_3, f)$, we denote it by $\{a,b,c\} \in \Delta$ and call $\{a,b,c\}$ forms
a triangle. For a spherical angle  structure $x: C(S, T) \to (0, \pi)$, we use $x_r$ to denote the value of $x$ at the
$r$-th corner and consider $x=(x_1, ...., x_n)$ as a vector in $\bold R^n$. Under this identification,  the space of
all spherical angle  structures
$AS(S, T) =\{ x \in (0, \pi)^n | $ whenever $r,s,t$ form a triangle, $(x_r, x_s, x_t) \in M_3$\} becomes
an open convex polyhedron of dimension $n$.  The capacity of the spherical angle  structure $x$, denoted by
$\Theta(x)$,  is  given by,
$$\Theta(x) = \sum_{ \{r,s,t\} \in \Delta} \theta(x_r, x_s, x_t).$$

Since $\theta(x_1, x_2, x_3)$ is strictly convex, we have,

\medskip
\noindent
{\bf Lemma 2.2.} \it The capacity function $\Theta$ defined on $AS(S, T)$ is a strictly convex function. \rm

\medskip
\noindent
2.5. Given any  map $D: E \to (0, \infty)$, we denote $AS(S, T;D)$ the subspace of all spherical angle  structures with
edge invariant  equal to $D$. 
\medskip
\noindent
{\bf Lemma 2.3.} \it If $AS(S,T;D) $ is non-empty,   then the critical points of $\Theta|_{AS(S, T;D)}$ are exactly those
spherical angle  structures derived from spherical polyhedron metrics. \rm
\medskip
\noindent
{\bf Proof. }  For simplicity, let us set $G = \Theta |_{AS(S, T;D)}$.  Applying the Lagrangian multipliers to $\Theta$ on $AS(S,T)$ subject
to the set of linear constraints $D_x(e) = D(e)$ for $e \in E$, we see that at a critical point of $G$, there is a 
map $C: E \to \bold R$ (the multipliers) so that, for all indices $i$,

$$ \partial \Theta/\partial x_i = C_e   \tag 2.4$$
where the i-th corner is of the form $(e,f)$, i.e., the i-th corner is facing the edge $e$.
Let the three corners of the triangle $f$ be labeled by $i,j,k$. Then $\partial \Theta/\partial x_i = \ln \tan ( y_i/2)$ where $y_i$
is given by the cosine law (2.1). This shows, by (2.4),  that the edge length of $e$ in the spherical triangle of inner angles
$x_i, x_j, x_k$ depends only on $C_e$. In particular, if $f'$ is the second triangle in $T$ having $e$ as an edge, then
the length of $e$ calculated in $f'$ in the spherical angle  structure is the same as the length of $e$ calculated using $f$.
In summary, we see that there is a well defined assignment of edge lengths $l: E \to (0, \pi)$ so that the assignment
on the three edges of each triangle forms the lengths  of a spherical triangle and the inner angles induced by $l$ is $x$.

To see the result in the other direction, suppose we have a point in $AS(S, T;D)$ which is induced from a spherical polyhedron metric 
$l : E \to (0, \infty)$. We want to show that the point is a critical point of $G$.
Since the constraints $D_x =D$ are linear, the critical points $p$ of $G$ on $AS(S, T;D)$ are the same as those points
$q \in AS(S, T;D)$ so that there is a map $C: E \to \bold R$ satisfying (2.4) at $q$. 
Evidently at a spherical angle  structure derived from a spherical polyhedron metric  $l: E \to \bold R_{>0}$, we define $C_e$ to be $\ln \tan ( l(e)/2)$. Then (2.4) follows. QED

\medskip
It is well known that for a smooth strictly convex function $f$ defined on a convex open $W$ set in $\bold R^n$, the gradient of $f$
is a diffeomorphim from $W$ to an open set in $\bold R^n$.
As a consequence of lemma 2.2 and lemma 2.3, we see theorem 1.1 follows. 

\medskip
\noindent
2.6. To prove theorem 1.2, by proposition 3.1, the function $\Theta$ on the space of all spherical angle  structure $AS(S, T;D)$
has a continuous extension to the closure $\bar AS(S, T;D)$ of $AS(S, T;D)$ in $\bold R^n$. The closure is evidently
compact since it is contained in $[0, \pi]^n$. Take a minimal point $a$ of $\Theta$ in the closure $\bar AS(S,T;D)$. If
the point $a$ is in $AS(S, T;D)$, we are done. We claim that  $a \in \partial AS(S,T;D)$ is impossible.
Suppose otherwise, there is a triple of indices $\{u',v', w'\}$
so that $(a_{u'}, a_{v'}, a_{w'})$ is in the boundary of $M_3$. Take a point $p \in AS(S,T;D)$ and consider the smooth path
$\gamma(t) = (1-t)a + t p$ for $t \in (0, 1] $ in $AS(S,T;D)$.  Let $g(t) = \Theta( \gamma (t))$. We have
$g(t) \geq g(0)$ for all $ t >0$ by the choice of the point $a$. Thus, $\liminf_{ t \to 0^+} dg/dt \geq 0$.
But, by proposition 4.1, we have $$\lim_{t \to 0^+} dg/dt = -\infty. \tag 2.5$$ This produces a contradiction.  
Here is the more detailed argument  to see (2.5).

Let $\Delta_1$ be the set of all triples of indices $\{ \{u,v,w\} | $ so that $\{u,v,w\} \in \Delta$ and $(a_u, a_v, a_w) \in \partial M_3$\}
and $\Delta_2 = \Delta -
\Delta_1$. Then the function $g$ can be written as,
$$ g(t) =\sum_{ \{u,v,w\} \in \Delta_1} \theta (x_u(t), x_v(t), x_w(t))   +
 \sum_{ \{u,v,w\} \in \Delta_2} \theta (x_u(t), x_v(t), x_w(t))$$
where $x(t) = x(\gamma (t))$. The derivative $g'(t)$ can be expressed as,
$$ g'(t) =\sum_{ \{u,v,w\} \in \Delta_1} d/dt [\theta (x_u(t), x_v(t), x_w(t))]   + \sum_{ \{u,v,w\} \in \Delta_2} d/dt[\theta (x_u(t), x_v(t), x_w(t))].$$
Note that since the edge invariant $D$ is assumed to be strictly less than $\pi$,  if  $\{u,v,w\}$ is in $\Delta_1$, then the triple 
  $(a_u, a_v, a_w)$ is in $\partial M_3 -\{ (0,0,\pi), (0,\pi, 0), (\pi,0,0), (\pi, \pi, \pi)\}$. Thus by proposition 
4.1, as $t$ tends to 0, each terms in the first sum tends to  $-\infty$. Each term in the second sum tends to a finite number as
$t$ tends to 0. Thus we see (2.5) holds.

\medskip
\noindent
2.7. The above proof in fact shows the following stronger result. \it A cycle \rm in the
triangulated surface (S,T) is an ordered collection of edges and triangles $\{ e_1,f_1, $$e_2,$$ f_2,  ..., e_n, f_n\}$ so that 
$e_i$ and $e_{i+1}$ are edges in $f_i$ and $e_1, e_n$ are edges of $f_n$.  An edge invariant assignment $D$
is said to contain a \{0,0,$\pi$\}-cycle if there is a cycle of edges and a point $a \in \partial AS(S, T; D)$ so that
$D_a(e_i) = \pi$ and the inner angles of each $f_i$ in $a$ are $0, 0, \pi$.

\medskip
\noindent
{\bf Theorem 2.1.} \it Given any triangulated surface and any edge invariant function $D: E \to (0, \pi]$ which contains no
$\{0,0,\pi\}$-cycles, if there is a linear
spherical structure having $D$ as the edge invariant, then there exists a spherical polyhedron metric  having $D$ as the edge invariant function. \rm
\medskip

The proof is evident.

\medskip
\noindent
\S 3. {\bf  Continuous Extension of the Capacity Function }
 
\medskip
We show that
the capacity of  spherical triangles extends  continuously to the degenerated triangles.  For the rest of the section, we take a spherical or hyperbolic triangle of inner angles
$x_1, x_2, x_3$ and edge lengths $y_1, y_2, y_3$ so that $y_i$-th edge is facing the $x_i$-th inner angle. We use $x=(x_1, x_2, x_3)$ and
 $x^*_i = 1/2( \pi + x_i - x_j -x_k)$.
As a convention, we assume the
indices $\{i,j,k\}=\{1,2,3\}$.
The main result of the section is the following.

\medskip
\noindent
{\bf Proposition 3.1.} \it  The capacity function $\theta(x) =
  \int_{(\pi/2, \pi/2, \pi/2)}^x \sum_{i=1}^3 \ln \tan (y_i/2) dx_i$ is given by the following,

$$ \theta (x_1, x_2, x_3) = - \sum_{i=1}^3 \Lambda (x^*_i) - \Lambda(  (\pi +x_1+ x_2 + x_3)/2)  + 4 \Lambda(\pi/4)\tag 3.1$$
and  the capacity function $\tilde \theta(x) =\int_{(0,0,0)}^x \sum_{i=1}^3 \ln \tanh(y_i/2) dx_i$ is given by
$$\tilde  \theta (x_1, x_2, x_3) = - \sum_{i=1}^3 \Lambda (x^*_i) - \Lambda((\pi+x_1+ x_2 + x_3)/2). \tag 3.2$$
In particular, both $\theta$ and $\tilde \theta$ have continuous extensions to the closure $\bar M_3$ of the moduli space of spherical triangles
$M_3  =\{(x_1, x_2, x_3) \in (0, \pi)^3 | x_1+x_2+x_3 > \pi$ and $x_i^* >0$, $i=1,2,3$\} and the closure of $\{(x_1, x_2, x_3)
\in (0, \pi)^3 | x_1+x_2+x_3 < \pi$\}. 
Geometrically, $16 \Lambda(\pi/4) -4\theta(x_1, x_2, x_3)$ is the volume of the hyperbolic ideal octahedron whose vertices are the intersection points of the three circles bounding the spherical triangle $(x_1, x_2, x_3)$.\rm

\medskip
\noindent
{\bf Proof. } The proof is a straight forward computation using the cosine law. Recall that the cosine law (2.1) says
$$ \cos y_i =\frac{ \cos  x_i + \cos  x_j \cos  x_k}{ \sin   x_j \sin   x_k }.$$
Use the summation formulas for cosine function that

$$ \cos (a+b) =\cos a \cos b - \sin a \sin b,$$
$$ \cos(a-b) = \cos a \cos b  +\sin a \sin b,$$
$$ \cos a + \cos b = 2 \cos ((a+b)/2) \cos((a-b)/2),$$
$$ \cos a - \cos b = 2 \sin ((a+b)/2) \sin ((b-a)/2),$$
we can rewrite the cosine law as one of the following,

$$ \cos y_i -1=  2 \frac{ \sin  x_i^* \cos( (x_i+x_j+x_k)/2)  }{\sin   x_j \sin   x_k}, \tag 3.3$$
and


$$ \cos y_i +1= 2 \frac{ \sin  x_j^* \sin x_k^*}{\sin   x_j \sin   x_k}.  \tag 3.4$$

In particular,
$$ \frac{1 -\cos y_i }{1 + \cos y_i } =  -\frac{  \sin  x_i^* \cos( (x_i+x_j+x_k)/2)}  { \sin  x_j^* \sin x_k^*}. \tag 3.5$$

However, we also have the trigonometric identity,
$$ \tan^2  (u/2) = \frac{ 1 -\cos u }{1 + \cos u}. $$
This shows that the cosine law for spherical triangles can be written as,
$$ \tan^2  (y_i/2)  = -\frac{  \sin  x_i^* \cos( (x_i+x_j+x_k)/2)}  { \sin  x_j^* \sin x_k^*}. \tag 3.6$$
By the same calculation and using $\tanh^2(u/2) = (\cosh u -1)/(\cosh u +1)$, we obtain the cosine law for hyperbolic triangles as,
$$ \tanh^2  (y_i/2)  = \frac{  \sin  x_i^* \cos( (x_i+x_j+x_k)/2)}  { \sin  x_j^* \sin x_k^*}. \tag 3.7$$

Since by definition, $\partial \theta/\partial x_i = \ln \tan (y_i/2)$, by (3.6), we have

$$ \partial \theta/\partial x_i = 1/2 [ \ln \sin x^*_i - \ln \sin x_j^* - \ln \sin x_k^*  + \ln (|\sin((x_1+x_2 + x_3+\pi)/2)|) ]\tag 3.8$$

Since the function $F(x_1, x_2, x_3)$ given by the right hand side of the (3.1) has the partial derivative,
$$ \partial F/\partial x_i = 1/2[  \ln (2 \sin x_i^*) - \ln (2 \sin x_j^*) - \ln (2 \sin x_k^*)  +\ln ( 2 |\sin ((x_1+x_2 + x_3+\pi)/2)|)]$$
$$=  1/2[ \ln \sin x^*_i - \ln \sin x_j^* - \ln \sin x_k^*  + \ln (|\sin((x_1+x_2 + x_3+\pi)/2)|)],$$
we see that $\partial F/\partial x_i = \partial \theta /\partial x_i$. In particular, these two functions differ by a constant
on $M_3$. Since $\theta(\pi/2, \pi/2, \pi/2) =0 =F(\pi/2, \pi/2, \pi/2)$, the result follows.  In particular, we see that
$\theta$ has a continuous extension to the 3-space $\bold R^3$. The same calculation using (3.7) verifies (3.2). 

Since  three great circles bounding a spherical triangle decompose the 2-sphere into eight spherical triangles, it follows that
the convex hull of the  six intersection points of three circles is the union of eight hyperbolic tetrahedra each of them has three vertices at the sphere at infinity and one vertex the Euclidean center.
By  (3.1) and known formula for volume of hyperbolic tetrahedra with three vertices at the sphere at infinity [Vi], i.e., (3.10) below,  we see that $16 \Lambda(\pi/4)-4 \theta(x_1, x_2, x_3)$ is the volume of the hyperbolic octahedron which is the convex hull of the six points. QED

\medskip
\noindent
{\bf Remarks 3.1.}
Proposition 3.1 shows that 
 the functions $\theta (x_1, x_2, x_3)$  and $\tilde \theta(x_1, x_2, x_3) $ are essentially $W(x)$ where,
$$ W(x_1, x_2, x_3)=  - \sum_{i=1}^3 \Lambda (x^*_i) - \Lambda((\pi+x_1+ x_2 + x_3)/2) \tag 3.9$$
This function $W(x)$ is closely related to 
$$V(x_1, x_2, x_3) = \sum_{i=1}^3( \Lambda(x_i)  + \Lambda (x^*_i)) -\Lambda(  (\pi +x_1+ x_2 + x_3)/2).  \tag 3.10$$
For a spherical triangle $x$, the function
$V(x)/2$  is known to be the hyperbolic volume of a hyperbolic tetrahedron with
 three vertices at the sphere at infinite so that the link at the finite vertex is the spherical $x$
(see [Vi], also [Le]). For  a hyperbolic triangle $(x_1, x_2, x_3)$, $V(x)$
is the volume of the convex hull of the intersection points of circles bounding the triangle. This is the
function used by Leibon as the capacity. For a Euclidean triangle x, $V(x)/2$ is the volume of the
hyperbolic ideal tetrahedron with dihedral angles $x_1, x_1, x_2, x_2, x_3, x_3$.  Peter Doyle [Le] noticed that $V(x)$ is not concave on $M_3$ and took $V(x)$ as a different capacity for spherical triangles.  He observed that the critical point of this capacity for spherical angle structures with prescribed Delaunay invariant are the spherical cone metrics. 
On the other hand,  $V(x)$ is concave in the set $\{ (x_1, x_2, x_3) \in [0, \pi]^3 | x_1 + x_2 + x_3 \leq \pi\}$ ([Le]). 
  For a spherical triangle $x$, $-4W(x)$ is the volume of
the ideal hyperbolic octahedron whose vertices are the intersection points of the circles bounding the triangle. For a Euclidean triangle $x$, we have $W(x)=-V(x)/2$. We do not know the geometric meaning
of $W(x)$ for  a hyperbolic triangle $x$.
The other related works are [CV] and [BS].
\medskip
\noindent
{\bf 3.2}. It can be shown that functions $W$ and $V$ in (3.9) and (3.10) are the only functions, up to scaling and adding of  linear functions, with the required properties. To be more precise, if $F(x_1, x_2, x_3)$ is a smooth function of the inner angles $(x_1, x_2, x_3)$ of a triangle so that
$\partial F/\partial x_i$ is a universal function of the edge length $y_i$, then $F = c_1W + c_2(x_1+x_2+x_3)+c_3$ for some constants 
$c_1, c_2$ and $c_3$. Similarly, if $F(x_1, x_2, x_3)$ is a smooth function so that $\partial F/\partial x_i^*$ is a 
universal function of $y_i$, then $F = c_1V+c_2(x_1+x_2+x_3) + c_3$ for some constants $c_1, c_2$ and $c_3$. This shows that if one intends to 
find the constant curvature cone metrics in $AS(S,T;D)$, $AS(S, T; \Cal D)$, $AH(S, T;D)$ or $AH(S, T; \Cal D)$ by a variational method so that the energy is contructed locally by summing up the energies of the triangles, then all the possible candidates of the energies are $c_1V+c_2(x_1+x_2+x_3)+c_3$ and $c_1W+c_2(x_1+x_2+x_3) + c_3$.

\medskip

\medskip
\noindent
\S 4. {\bf Degeneration of Spherical Triangles}
\medskip
\noindent
The goal of this section is to understand how a sequence of spherical triangles degenerates and to understand the behavior of the
derivatives of the capacity on the sequence of degenerated spherical triangles. 
Recall that the moduli space $M_3$ of spherical triangles is an open regular tetrahedron in the 3-space.
The closure $\bar M_3$ of $M_3$ is the closed tetrahedron. We call a point in the boundary $\partial M_3 = \bar M_3 - M_3$
a \it degenerated \rm spherical triangle (with respect
to inner angles).  The goal of the section is to prove,
\medskip
\noindent
{\bf Proposition 4.1.} \it 
For any  point $a \in \bar M_3 - M_3$ and a point $p \in M_3$, let  $f(t) =\theta ( (1-t) a + t p)$
where $t \in [0,1]$. If $a $ is not equal to any of the points $ (0,0,\pi), (0, \pi, 0), $$(\pi, 0, 0),$$ (\pi, \pi, \pi)$, then
$$\lim_{t \to 0^+} f'(t) = -\infty. \tag 4.1$$
If $a \in \{ (0,0,\pi), (0, \pi, 0), (\pi, 0, 0), (\pi, \pi, \pi)\}$, then the limit $\lim_{t \to 0^+} f'(t)$ exists and is a finite number. 
 \rm


\medskip
\noindent
4.1.  The moduli space $M_3$ of spherical triangles is given by
 $\{ x \in (0, \pi)^3 | x_i^* >0, x_1+x_2+x_3 > \pi\}$ which is the open regular tetrahedron inscribed in the standard cube
$[0, \pi]^3$. 
The four vertices of the tetrahedron are $v_1=(\pi,0,0), v_2 =(0,\pi,0), v_3=(0,0,\pi)$ and $v_4 =(\pi,\pi, \pi)$
and its four triangular faces lie in the planes given by the
linear equations $x^*_i=0$, i=1,2,3, and $x_1+x_2+x_3 =\pi$ respectively. We now decompose the boundary $\partial M_3$ into a disjoint union 
of six parts, denoted by $I,II,III,IV,V$ and $VI$, as follows. Here $I$ is the open triangle $\Delta v_1 v_2v_3$. Part $II$ is the union of the
three open triangles $\Delta v_4v_iv_j$ where $\{i,j\} \subset \{1,2,3\}$. Part $III$ is the union of three open edges of the triangle
$I$, i.e., $III$ is the union of open intervals $v_iv_j$ where $\{i,j\} \subset \{1,2,3\}$ . Part $IV$ is the union of the three open intervals
$v_4v_i$. Part $V$ is \{$ (\pi,\pi,\pi)$\}. Part $VI$ is $\{ (0,0,\pi), (0,\pi,0), (\pi,0,0)\}$. The algebraic description of them is as follows.

$$ I =\{ a \in (0, \pi)^3 | a_1 + a_2 + a_3 =\pi, a^*_i \in (0, \pi)\},$$
$$ II = \cup_{i=1}^3 \{ a \in (0, \pi)^3 |  a_i^* =0, a^*_j, a_k^* \in (0, \pi),   a_1+a_2+a_3 >\pi \},$$
$$ III =\cup_{i=1}^3 \{ a \in [0, \pi)^3 |  a^*_i=0, a_j^*, a_k^* \in (0, \pi),  a_1+a_2+a_3 =\pi \},$$
$$ IV =\cup_{i=1}^3 \{ a \in (0, \pi]^3 | a_i =\pi,  a^*_j =a^*_k =0, a^*_i \in (0, \pi),  a_1 + a_2+a_3 >\pi\},$$
$$ V =\{ a \in [0, \pi]^3 | a_i^* = 0, i=1,2,3, a_1+a_2+a_3 = 3\pi \},$$
$$ VI=\cup_{i=1}^3 | a^*_j = a^*_k =0, a^*_i=2\pi,  a_1+a_2+a_3 =\pi \}.$$

\medskip
As usual, we have used the convention that $\{i,j,k\}=\{1,2,3\}$ above.

\medskip
\noindent
4.2. We now prove proposition 4.1 by considering the limit $\lim_{t \to 0^+} f'(t)$ according to the type of the degenerated spherical
triangle $a$.  Let $a=(a_1, a_2, a_3), p=(p_1, p_2, p_3)$ and let  $x_i = x_i(t) = (1-t) a_i + t p_i$. We use $y_i = y_i(t)$
to denote the corresponding edge lengths of the triangle $x = (x_1, x_2, x_3)$.  Note that, $x_i \to a_i$ and $x_i^* \to a_i^*$
as time $t$ tends to 0, also $dx_i/dt = p_i - a_i$.  By definition,
$$ f'(t) =  \sum_{i=1}^3 \ln \tan (y_i(t)/2) (p_i -a_i).  \tag 4.2$$
By (3.8), we write,
$$\ln \tan (y_i/2) =   S(x^*_i) - S(x^*_j)  -S(x^*_k) + C(x) \tag 4.3$$
where $S(u) = 1/2 \ln \sin (u)$ and $C(x) = 1/2 \ln |\cos((x_1+x_2+x_3)/2 )|$.  Assume in the following computation that
$(i,j,k)$ is a cyclic permutation of $(1,2,3)$, or more precisely, we take $j=i+1, k=i+2$ where indices are counted modulo 3.
Substitute (4.3) into (4.2), we obtain,
$$ f'(t) = \sum_{i=1}^3 (  S(x^*_i) - S(x^*_{i+1})  -S(x^*_{i+2}) + C(x) ) (p_i -a_i)$$
$$= \sum_{i=1}^3   (S(x^*_i) - S(x^*_{i+1})  -S(x^*_{i+2}) ) (p_i- a_i) + C(x)  (\sum_{i=1}^3p_i - \sum_{i=1}^3 a_i)$$
$$= 2 \sum_{i=1}^3  S(x^*_i) ( p^*_i - a^*_i) + C(x)  (\sum_{i=1}^3 p_i - \sum_{i=1}^3 a_i) \tag 4.4$$

We  now discuss the limit of $f'(t)$ as $t$ tends to 0 according to the type of the degenerated triangle $a$.

\medskip
\noindent
4.3. Case 1, the triangle $a$ has type I, i.e., $a_1+a_2+a_3 =\pi$ and $a_i, a_i^* \in (0, \pi)$.
In particular, $\lim_{t \to 0^+}S( x^*_i ) = S(a_i^*)$ exists in $\bold R$.  Thus the unbounded term in (4.4) is
the last term $C(x)  (\sum_{i=1}^3 p_i - \sum_{i=1}^3 a_i) $ which tends to $-\infty$ due to
 $a_1+a_2+ a_3 = \pi$,  $p_1+p_2+p_3 > \pi$ and $\lim_{t \to 0^+} C(x) =-\infty$. This shows the proposition for case 1.

\medskip
\noindent
4.4. Case 2, the triangle $a$ has type II. For simplicity, we may assume that $ \pi + a_1 = a_2 + a_3$, i.e., $a^*_1=0$, 
$a_i, a_2^*, a_3^* \in (0, \pi)$, and $a_1+a_2+a_3 \in (\pi, 3\pi)$.  Then  the unbounded term in (4.4) is 
$ 2S(x^*_1) (p_1^* - a^*_1)$. All other terms are bounded since the $\lim_{t \to 0^+}S( x^*_i(t)) = S( a^*_i)$
is finite for $i=2,3$ and $\lim_{t \to 0^+} C(x) = 1/2 \ln |\cos (a_1+a_2+a_3)/2)|$ is also finite. On the other hand, $p_1^* >0$
, $a^*_1=0$ and $\lim_{t \to 0^+} S(x^*_1) = -\infty$, we see that $\lim_{t \to 0^+} f'(t) = -\infty$.

\medskip
\noindent
4.5. Cases 3,4, the triangle $a$ has type III or IV. In these cases, exactly two of the four equations $a^*_1=0, a_2^*=0,
a_3^*=0, $ or $a_1+a_2+a_3 =\pi$ hold.  To be more precise, in the case III, we may assume without loss of generality that
$a_1^*=0$, $\sum_{i=1}^3 a_i = \pi$, $a_2^*, a_3^* \in (0, \pi)$. Thus, in (4.4), exactly two terms, $2S(x_1^*)(p_1^*-a_1^*)$
and $C(x)(\sum_{i=1}^3 p_i -\pi)$ tend to $-\infty$ as $t$ approaches 0. The other two terms remain bounded. Thus the result
follows.

In the case IV, we may assume for simplicity that $a_1^*=a_2^*=0$ and $\sum_{i=1}^3 a_i >\pi$ and $a^*_3 >0$. Then
due to $0 < \sum_{i=1}^3 a_i^* = 3\pi -\sum_{i=1}^3 a_i$, we have $\sum_{i=1}^3 a_i < 3\pi$. This shows that
$\lim_{t \to 0} C(x) =C(a)$ is finite. Thus in (4.4), there are again exactly two terms, namely $2 S(x_1^*)(p_1^*-a_1^*)$
and $2 S(x_2^*)(p_2^*-a_2^*)$ tend to $-\infty$ as $t$ approaches 0. The other two terms remain bounded. Thus the result follows
again.

\medskip
\noindent
4.6. Case 5, the triangle $a$ is an equator $(\pi, \pi, \pi)$. In this case $a^*_i =0$ and $a_1+a_2+a_3 = 3\pi$. 
Using (4.2) and (4.3), we have,
$$ f'(t) = \sum_{i=1}^3 (S(x^*_i) -S(x_j^*) - S(x_k^*)   + C(x) ) (p_i- \pi)$$
$$ = \sum_{i=1}^3 [(S(x^*_i)-S(x_j^*)) + ( C(x) - S(x^*_k))] (p_i - \pi) \tag 4.5$$

We note that both limits $\lim_{t \to 0^+}( S(x_i^*) -S(x_j^*))$ and $\lim_{t \to 0^+}( C(x) -S(x_k^*)) $ exist in $\bold R$. 
Indeed, by definition,

$$ x^*_i =1/2[ \pi + x_i - x_j -x_k] =1/2[ \pi + (1-t) ( a_i - a_j -a_k) + t ( p_i -p_j -p_k)]$$
$$=1/2[ \pi + (1-t) (-\pi) + t(p_i - p_j -p_k)$$
$$= 1/2( t ( p_i -p_j -p_k + \pi)) = t p_i^*.$$
$$ x_1+x_2+x_3 = t(p_1+p_2 + p_3) + (1-t) 3\pi = 3\pi + t( p_1+p_2 + p_3 - 3\pi)$$

Thus, $ S(x^*_i) -S(x_j^*) = 1/2( \ln \sin(  tp_i^*)  - \ln \sin (t p_j^*))$ which tends to $ 1/2(\ln \sin p_i^* - \ln \sin p_j^*)$ as $t$ tends to 0. Similarly, $ C(x) - S(x^*_k) $ tends to  the finite number $1/2(\ln |\sin((p_1+p_2+p_3 - 3\pi)/2)|  - \ln ( \sin( p_k^*)))$.

\medskip
\noindent
4.7. Case 6, the triangle $a$ is of type VI. For simplicity, we assume that $a=(\pi, 0, 0)$. Thus $a_1+a_2+a_3=\pi$,
$a_1^* = 2\pi$, $a_2^*=a_3^*=0$.  We use (4.5) to calculate the limit $\lim_{t \to 0} f'(t)$. The calculation is exactly the same as
that of case 5. 
Indeed,  each of the four terms $S(x^*_i)$ and $C(x)$ tends to $-\infty$ as $t$ approaches zero. On the other hand, by the
same argument as in 4.6,  both of 
the limits 
$ \lim_{t \to 0^+} S(x^*_i)/ S(x_k^*)$ and $\lim_{t \to 0^+} S(x^*_i)/C(x)$ are finite. Thus the result follows.

This ends the proof of proposition 4.1.
\medskip
\noindent


\medskip
\noindent
4.8. {\bf Remark.} We give a geometric interpretation of the stratification I, II, ..., VI of the degenerated triangles.
The type I boundary point $x \in $ $\{ x \in (0, \pi)^3 | x_1 + x_2 + x_3 = \pi$\} corresponds to the "Euclidean triangle".
Geometrically, it represents a point which is the limit of  spherical triangles  shrinking to a point
 so that its inner angles tend to three numbers in $(0, \pi)$. In particular, if one defines the edge length
$y_i=0$ for these triangle, the cosine law (2.1) still makes sense in terms of taking limit. 
The type II points in $\{ x \in (0, \pi)^3 | x_1+x_2 + x_3 > \pi,$ $x^*_i=0$, $x^*_j >0$, $x^*_k >0\}$ correspond
 to the other codimension-1 faces. They represent  the "exceptional Euclidean triangles". 
 Geometrically,
it  is  the limit of sequence of spherical triangles expanding to a union of two geodesics from a point to its antipodal 
point so that the inner angles tend to three numbers in $(0, \pi)$. 
In particular, the edge lengths are $y_i=0, y_j = y_k=\pi$ and a type II triangle
has two vertices.
Note that the edge length function $y_i$ extends continuously on the set $M_3 \cup I \cup II$.
There are two types of codimension-2 faces. The first type, denoted by III, consists of three open edges of the
form $\{ x =(x_1, x_2, x_k) \in [0, \pi)^3 | $ $x_i=0$, $x_j, x_k >0$ and $x_j+x_k=\pi$\}. This is a further degeneration
of "Euclidean triangles". The second type of codimension-2 face, denoted by IV, consists of the three open edges of
the form $\{ x =(x_1, x_2, x_3 ) \in (0, \pi]^3 | x_i = \pi$, $x_j=x_k \in (0, \pi)$\}.  Geometrically, it corresponds to
a degenerated spherical triangle so that two of its three distinct vertices are antipodal points. Due to the location of
the third vertex (of inner angle $\pi$), the length functions $y_r$ does not extend continuously from $M_3$ to $M_3 \cup IV$.
Finally, there are two types of vertices. The first type, denoted by $V$,  is the point $(\pi, \pi, \pi)$
 corresponding to the equator and the second type, denoted by $VI$, 
consists of $(0, 0, \pi), (0, \pi, 0), (\pi, 0, 0)$ corresponding to a degenerated triangle whose 
three distinct  vertices lie in a great circular arc of length at most $\pi$.

\medskip
\noindent
\centerline{\bf References}

[BS]
Bobenko, Alexander I.; Springborn, Boris A. Variational principles for circle patterns and Koebe's theorem. Trans. Amer. Math. Soc. 356 (2004), no. 2, 659--689, MR2121737, Zbl 1044.52009.

[CV] Colin de Verdiere, Yves, Un principe variationnel pour les empilements de cercles. Invent. Math. 104 (1991), no. 3, 655--669, MR1106755,
Zbl 0745.52010.

[Gu] Guo, Ren,  Geometric angle structures on triangulated surfaces, 
\newline http://front.math.ucdavis.edu/math.GT/0601486.

[Le] Leibon, Gregory, Characterizing the Delaunay decompositions of compact hyperbolic surfaces. Geom. Topol. 6 (2002), 361--391,
MR1914573, Zbl 1028.52014.

[Lu1] Luo, Feng,  Volume and angle structures on 3-manifolds, preprint, 2005,
\newline http: front.math.ucdavis.edu/math.GT/0504049.

[Lu2] Luo, Feng,  Continuity of the volume of simplices in classical geometry,  to appear in Comm. Cont. Math.  http://front.math.ucdavis.edu/math.GT/0412208.

[Mi] Milnor, John, computation of volume, chapter 7 of Thurston's note on geometry and topology of 3-manifolds, 
at www.msri.org/publications/books/gt3m/.

[Ri1] Rivin, Igor, Euclidean structures on simplicial surfaces and hyperbolic volume. Ann. of Math. (2) 139 (1994), no. 3, 553--580.
MR1283870, Zbl 0823.52009.

[Ri2] Rivin, Igor,  Continuity of volumes -- on a generalization of a conjecture of J. W. Milnor.  http://front.math.ucdavis.edu/math.GT/0502543.

[Th] Thurston, William, the geometry and topology of 3-manifolds, lecture notes, Math Dept., Princeton University,   1978,                                       
 at www.msri.org/publications/books/gt3m/.

[Vi]  Vinberg, E. B., The volume of polyhedra on a sphere and in Lobachevsky space. Algebra and analysis (Kemerovo, 1988), 15--27, Amer. Math. Soc. Transl. Ser. 2, 148, Amer. Math. Soc., Providence, RI, 1991, MR1109060, Zbl 0742.51019.

\medskip
\noindent

Department of Mathematics

Rutgers University

Piscataway, NJ 08854, USA

email: fluo\@math.rutgers.edu

\end

\end